\documentclass[12pt,a4paper]{article}
\usepackage{fullpage}
\usepackage{amssymb, amsmath, amsfonts, amsthm}
\usepackage[english]{babel}
\usepackage{indentfirst}
\usepackage[parfill]{parskip}
\usepackage{hyperref, color}

\author{\large\textbf{Vadim Demichev}\\ \normalsize\emph{The University of Sheffield}\\ \normalsize \emph{S10 2TN Sheffield, UK}\\ \normalsize \emph{vadim.demichev@gmail.com}}
\title{\vspace{-30pt}\large\textbf{On the invariance principle\\ for empirical processes of associated sequences}}
\date{}

\begin{document}
\maketitle

\begin{abstract}
We consider empirical processes generated by strictly stationary sequences of associated random variables. S.\,Louhichi established an invariance principle for such processes, assuming that the covariance function decays rapidly enough.  We show that under certain conditions imposed on the pairwise distributions of the random variables in question the restrictions on the rate of decay of the covariance function can be relaxed. \par
\end{abstract}

\vspace{0pt}\textbf{Introduction}\vspace{0pt}

The theory of associated random processes and fields has been actively developed since 1980s. The property of association is easy to verify for a wide range of random fields arising in a number of applications (see e.g. \cite{BS2007}). Furthermore, limit theorems under association can often be formulated with conditions imposed only on the covariance function, and the latter is usually not difficult to estimate (see \cite{BS2007}, \cite{O2012}, and \cite{P2012} for an extensive review of asymptotic results for associated random variables). It is, therefore, often convenient to carry out asymptotic analysis of random processes and fields via the use of the property of association. Thus, limit properties of associated sequences are of considerable interest. \par

The classical invariance principle in the Skorokhod space for empirical distribution functions was generalised to the case of associated random variables in \cite{Y1993}. The restrictions on the rate of decay of the covariance function were later weakened first in \cite{SY1996} and then in \cite{L2000}. We show that these restrictions can be further relaxed, provided the two-dimensional distributions exhibit a certain type of positive dependence stronger than association. Note that related invariance principles in the spaces ${\sf{L}}_p$, $p>2$, and ${\sf{L}}_2$ were obtained in \cite{OS1998} and \cite{MS2002}, respectively.

\vspace{0pt}\textbf{Invariance principle}\vspace{0pt}

Let us recall the necessary definitions. Consider a family of random variables ${\sf{X}}=\{X_t,\;t\in T\}$. According to \cite{EPW1967} ${\sf{X}}$ is called {\it associated} (we write ${\sf{X}}\in{\sf{A}}$) if for any finite $I,J\subset T$ and bounded coordinate-wise nondecreasing functions $f:\mathbb{R}^{card(I)}\rightarrow\mathbb{R}$, $g:\mathbb{R}^{card(J)}\rightarrow\mathbb{R}$ one has
$${\sf{cov}}(f(X_i,\,i\in I),\;g(X_j,\,j\in J))\geqslant 0.$$
Associated families of random variables possess an important property, which is very useful in establishing various moment estimates. For square-integrable ${{\sf{X}}=\{X_t,\;t\in T\}\in{\sf{A}}}$, any finite $I,J\subset T$, and arbitrary Lipschitz functions $f:\mathbb{R}^{card(I)}\rightarrow\mathbb{R}$, $g:\mathbb{R}^{card(J)}\rightarrow\mathbb{R}$ the following inequality holds (see e.g. \cite[Theorem 1.5.3]{BS2007})
\begin{equation}
{\sf{cov}}(f(X_i,\,i\in I),\;g(X_j,\,j\in J))\leqslant \sum_{i\in I}\sum_{j\in J}{\sf{Lip}}_i(f){\sf{Lip}}_j(g)|{\sf{cov}}(X_i,X_j)|, \label{a15}
\end{equation}
where ${\sf{Lip}}_l(f)$ is the Lipschitz constant with respect to $x_l$, $l\in I$, of a Lipschitz function $f(x_i,\,i\in I)$, $x_i\in\mathbb{R}$, $i\in I$.\par

Let $X,Y$ be real-valued random variables. $X$ is said to be {\it stochastically increasing} in $Y$ \cite{T1958} (this type of dependence is also called {\it positive regression dependence}) if ${{\sf{P}}(X>x|Y=y)}$ is a nondecreasing function of $y\in\mathbb{R}$ for all $x\in\mathbb{R}$. Following \cite{N2006}, we write ${\sf{SI}}(X|Y)$. \par

It is known (see e.g. \cite[Corollary 5.2.11]{N2006}) that if $X,Y$ are continuous, then ${\sf{SI}}(X|Y)$ if and only if the copula $C(u,v)={\sf{P}}(F_X(X)\leqslant u,F_Y(Y)\leqslant v)$ is a concave function of $v\in[0,1]$ for any $u\in[0,1]$. Here $F_X$ and $F_Y$ are the distribution functions of $X$ and $Y$, respectively. In addition, if the random vector $(X,Y)$ belongs to the ${\sf{TP}_2}$ class (sometimes also called {\sf{PLR}} class), then ${\sf{SI}}(X|Y)$ and ${\sf{SI}}(Y|X)$ (see e.g. \cite[Theorem~5.2.19]{N2006}). In particular, if $(X,Y)$ has a positive density $f\in C^2(\mathbb{R}^2)$, $\partial^2({\sf{log}}\;f)/\partial x\partial y\geqslant 0$, $x,y\in\mathbb{R}$, then ${\sf{SI}}(X|Y)$ and ${\sf{SI}}(Y|X)$ \cite[Lemma 1.4.6]{BS2007}. For example, any positively correlated Gaussian vector $(X,Y)$ is ${\sf{TP}_2}$ and hence possesses the latter property. \par

Let ${\sf{X}}=\{X_k,\;k\in\mathbb{Z}\}$ be a strictly stationary sequence of uniformly distributed on $[0,1]$ random variables, ${\sf{X}}\in {\sf{A}}$. Suppose that for some $C,\alpha>0$
\begin{equation}
{\sf{cov}}(X_0,X_k)\leqslant Ck^{-\alpha},\;\;\;k\in\mathbb{N}. \label{a1}
\end{equation}
Introduce the condition \\
{\sf{(SI)}} ${\sf{SI}}(X_k|X_l)$ and ${\sf{SI}}(X_l|X_k)$ for any $k\neq l\in\mathbb{Z}$. \par

Put
\begin{equation}
G_n(x)=\frac{1}{\sqrt{n}}\sum_{i=1}^n\left(\mathbb{I}\{X_i\leqslant x\}-{\sf{P}}(X_i\leqslant x)\right),\;\;\;x\in\mathbb{R}. \label{a6}
\end{equation}
Now we can formulate our main result. \par

\textbf{Theorem 1}. \emph{Assume {\sf{(SI)}}. Suppose also that $(\ref{a1})$ holds for some $\alpha>(5+\sqrt{17})/4\approx 2.28$. Then the random processes $G_n$ converge weakly in the Skorokhod space $D([0,1])$ to a centered Gaussian process $G$ with covariance function
\begin{equation}
{\sf{cov}}(G(x),G(y))=\sum_{k\in\mathbb{Z}}{\sf{cov}}(\mathbb{I}\{X_0\leqslant x\},\,\mathbb{I}\{X_k\leqslant y\}),\;\;\;x,y\in\mathbb{R}. \label{a7}
\end{equation}
}In \cite{L2000} {\sf{(SI)}} is not required, but the functional convergence of $G_n$ to $G$ is obtained under the more restrictive assumption $\alpha>4$. \par

Let $D([-\infty,\infty])$ \cite{DD2011} be the image of $D([0,1])$ under the bijection $f\longmapsto f\circ\phi$, $f\in D([0,1])$, where $\phi(x)=1/2+arctan(x)/\pi$, $x\in[-\infty,\infty]$. As in \cite{L2000} we can formulate the invariance principle for the case of an arbitrary continuous distribution. \par

\textbf{Corollary 1}. \emph{Let ${\sf{X}}=\{X_k,\;k\in\mathbb{Z}\}$ be a strictly stationary sequence of random variables with continuous distribution function $F$, ${\sf{X}}\in {\sf{A}}$. Assume {\sf{(SI)}}. Suppose also that ${\sf{cov}}(F(X_0),F(X_k))=\mathcal{O}(k^{-\alpha})$, $k\rightarrow\infty$, for some $\alpha>(5+\sqrt{17})/4$. Then the random processes $G_n$ $($defined by $(\ref{a6}))$ converge weakly in $D([-\infty,\infty])$ to a centered Gaussian process $G$ with covariance function of the form $(\ref{a7})$.} \par

Note that if $X_0$ is square-integrable, and $F$ is not only continuous but also Lipschitz, i.e. $X_0$ has a density bounded by certain $a>0$, then ${\sf{cov}}(F(X_0),F(X_k))\leqslant a^2{\sf{cov}}(X_0,X_k)$, $k\in\mathbb{N}$. This inequality follows from (\ref{a15}). \par

\textbf{Proof}. By \cite[Corollary 5.2.11]{N2006} we have ${\sf{SI}}(F(X_k)|F(X_l))$ and ${\sf{SI}}(F(X_l)|F(X_k))$, ${k\neq l\in\mathbb{Z}}$.
Therefore, the statement of Corollary 1 follows from Theorem 1 and the standard argument used to deduce the convergence in the general case from the convergence in the case of uniformly distributed random variables (see e.g. \cite[Theorem~16.4]{B1968}).
$\square$ \par

Let ${\sf{Lip}}(f)$ denote the Lipschitz constant of a Lipschitz function $f$ on $\mathbb{R}$. For $x,y\in\mathbb{R}$ put $x\vee y=\max\{x,y\}$. The proof of Theorem 1 is based on the following moment estimate. \par

\textbf{Lemma 1}. \emph{Assume {\sf{(SI)}} and suppose that $(\ref{a1})$ holds for some $\alpha>1$. Then for any $p>2$ and $\nu>0$ one can find such $K=K(p,\nu,\alpha,C)$ that
$${\sf{E}}\left|G_n(t)-G_n(s)\right|^p$$
\begin{equation}
\leqslant K\left(n^{(p/2-\alpha)\vee (1+\nu-p/2)}+|t-s|^{(1-1/\alpha)p/2}\right),\;\;\;s,t\in\mathbb{R},\;\;\;n\in\mathbb{N}. \label{a5}
\end{equation}
} \par

\textbf{Proof}. Fix $s<t\in\mathbb{R}$ and introduce the function
$$g(x)=\mathbb{I}\{x\leqslant t\}-\mathbb{I}\{x\leqslant s\}=\mathbb{I}\{x\in(s,t]\},\;\;\;x\in\mathbb{R}.$$
Clearly,
\begin{equation}
G_n(t)-G_n(s)=\frac{1}{\sqrt{n}}\sum_{k=1}^n(g(X_k)-{\sf{E}}g(X_k)),\;\;\;n\in\mathbb{N}. \label{a13}
\end{equation}
It is not difficult to show that for any finite disjoint $I,J\subset \mathbb{Z}$ the relation (\ref{a15}) implies
$${\sf{cov}}\left(\left|\sum_{i\in I}(g(X_i)-{\sf{E}}g(X_i))\right|,\left|\sum_{j\in J}(g(X_j)-{\sf{E}}g(X_j))\right|^{p-1}\right)$$
\begin{equation*}
\leqslant (p-1)|J|^{p-2}\sum_{u_1,u_2\in\{s,t\}}\sum_{i\in I}\sum_{j\in J}|{\sf{cov}}(\mathbb{I}\{X_i\leqslant u_1\},\mathbb{I}\{X_j\leqslant u_2\})|.
\end{equation*}
Here we used the inequality $\left|\sum_{j\in J}(g(X_j)-{\sf{E}}g(X_j))\right|\leqslant |J|$. Applying \cite[Theorem 1]{M2004} and (\ref{a1}), we get
\begin{equation}
|{\sf{cov}}(\mathbb{I}\{X_0\leqslant u_1\},\mathbb{I}\{X_k\leqslant u_2\})|\leqslant 4\,{\sf{cov}}(X_0,X_k)\leqslant 4\,Ck^{-\alpha},\;\;\;u_1,u_2\in\mathbb{R},\;\;\;k\in\mathbb{N}. \label{a3}
\end{equation} \par
Set $m={\sf{dist}}(I,J)=\min_{i\in I,j\in J}|i-j|$. In view of (\ref{a3})
$${\sf{cov}}\left(\left|\sum_{i\in I}(g(X_i)-{\sf{E}}g(X_i))\right|,\left|\sum_{j\in J}(g(X_j)-{\sf{E}}g(X_j))\right|^{p-1}\right)$$
\begin{equation}
\leqslant A\cdot(p-1)|I||J|^{p-2}Cm^{1-\alpha}, \label{a4}
\end{equation}
where $A=A(\alpha)>0$ depends only on $\alpha$.
Using (\ref{a4}) and directly following the proof of \cite[Lemma 1]{D2015} (cf. \cite[Theorem 4.2]{SY1996}), one can obtain the estimate
$${\sf{E}}\left|\sum_{k=1}^n(g(X_k)-{\sf{E}}g(X_k))\right|^p$$
\begin{equation}
\leqslant K_1\left(n^{1+\nu}+Cn^{(p-\alpha)\vee (1+\nu)}+n^{p/2}\left(\sum_{k\in\mathbb{Z}}|{\sf{cov}}(g(X_0),g(X_k))|\right)^{p/2}\right), \label{a11}
\end{equation}
where $K_1=K_1(p,\nu,\alpha)>0$ depends only on $p$, $\nu$, and $\alpha$.
\par

To show (\ref{a5}), set $\delta=|t-s|$ and note that $|{\sf{cov}}(g(X_0),g(X_k))|\leqslant \delta$, $k\in\mathbb{Z}$.
Therefore, by (\ref{a3})
$$\sum_{k=0}^{\infty}|{\sf{cov}}(g(X_0),g(X_k))|\leqslant \sum_{k\leqslant \delta^{-1/\alpha}}\delta+\sum_{k>\delta^{-1/\alpha}}16\,Ck^{-\alpha}\leqslant A_1\,\delta^{1-1/\alpha},$$
where $A_1=A_1(C,\alpha)>0$ depends only on $C$ and $\alpha$.
The latter inequality, (\ref{a11}), and (\ref{a13}) now yield (\ref{a5}).
$\square$ \par

\textbf{Proof of Theorem 1}. The convergence of finite dimensional distributions of $G_n$, $n\in\mathbb{N}$, follows from \cite[Theorem 4]{D2014}. It remains to show that the distributions of $G_n$, $n\in\mathbb{N}$, in $D([0,1])$ are tight. By \cite[Theorem 15.5]{B1968} it is sufficient to verify the relation
\begin{equation*}
\limsup_{n\rightarrow\infty}\;{\sf{P}}\left(\sup_{s,t\in[0,1],\;0<t-s<2^{-d}}|G_n(t)-G_n(s)|>\varepsilon\right)\xrightarrow[d\rightarrow \infty]{}0,\;\;\;\varepsilon>0.
\end{equation*}
For $m\in\mathbb{N}$, $t\in[0,1]$ put
$$t_m^-=\max\{x\in 2^{-m}\mathbb{Z}:\;x\leqslant t\},\;\;\;t_m^+=\min\{x\in 2^{-m}\mathbb{Z}:\;x\geqslant t\}.$$
It is easy to see that for any $0\leqslant s<t\leqslant 1$, $m\in\mathbb{N}$
$$G_n(t_m^-)-G_n(s_m^+)-2^{-m}\cdot 2\sqrt{n}$$
$$\leqslant G_n(t)-G_n(s)$$
$$\leqslant G_n(t_m^+)-G_n(s_m^-)+2^{-m}\cdot 2\sqrt{n},\;\;\;n\in\mathbb{N}.$$
Therefore, it is sufficient to show that for some nondecreasing sequence $m_n\in\mathbb{N}$, $n\in\mathbb{N}$, such that $2^{-m_n}=o(n^{-1/2})$, $n\rightarrow\infty$, the following relation holds
\begin{equation}
\limsup_{n\rightarrow\infty}\;{\sf{P}}\left(\sup_{s,t\in 2^{-m_n}\mathbb{Z},\;0<t-s<2^{-d}}|G_n(t)-G_n(s)|>\varepsilon\right)\xrightarrow[d\rightarrow \infty]{}0,\;\;\;\varepsilon>0.\label{a8}
\end{equation}
To estimate the probability in the left-hand side of (\ref{a8}), we use the chaining argument. Set
$$M_k=M_k(n)=\sup_{t\in 2^{-k}\mathbb{Z}}|G_n(t)-G_n(t-2^{-k})|,\;\;\;k\in\mathbb{N}.$$
It is easy to see that for $\mathbb{N}\ni d<m_n$ and arbitrary $s,t\in 2^{-m_n}\mathbb{Z}$, $0<t-s<2^{-d}$, we have
$$|G_n(t)-G_n(s)|\leqslant |G_n(t)-G_n(t_d^-)|+|G_n(s)-G_n(t_d^-)|\leqslant 2\sum_{k=1}^{m_n-d}M_{d+k},\;\;\;n\in\mathbb{N},$$
and, therefore,
$${\sf{P}}\left(\sup_{s,t\in 2^{-m_n}\mathbb{Z},\;0<t-s<2^{-d}}|G_n(t)-G_n(s)|>\varepsilon\right)\leqslant \sum_{k=1}^{m_n-d}{\sf{P}}\left(M_{d+k}>\varepsilon r^{k-1}(1-r)/2\right)$$
\begin{equation}
\leqslant \sum_{k=1}^{m_n-d}2^{d+k}\sup_{t\in 2^{-d-k}\mathbb{Z}}{\sf{P}}\left(|G_n(t)-G_n(t-2^{-d-k})|>\varepsilon r^{k-1}(1-r)/2\right),\;\;\;r\in(0,1). \label{a9}
\end{equation}
By Lemma 1 the right-hand side of (\ref{a9}) is not greater than
\begin{equation}
\sum_{k=1}^{m_n-d}2^{d+k}K\,2^p\frac{n^{(p/2-\alpha)\vee (1+\nu-p/2)}+\left(2^{-d-k}\right)^{(1-1/\alpha)p/2}}{r^{p(k-1)}\varepsilon^p(1-r)^p},\;\;\;p>2,\;\;\;\nu>0. \label{a10}
\end{equation}
Choose any $p\in\left((\alpha+1)\vee (2\alpha/(\alpha-1)),\;2\alpha-1\right)$. This interval is nonempty, since ${\alpha>(5+\sqrt{17})/4}$. Clearly, $p/2-\alpha>1-p/2$, and $\alpha>(p+1)/2$.
We can also choose $r=r(p,\nu)\in(0,1)$ large enough to satisfy $r^{-p}\leqslant 2^{\nu}$, $\nu>0$.
For $m_n=m_n(\nu)=\lfloor\log_{2}n^{1/2+\nu}\rfloor\vee 1$, $n\in\mathbb{N}$, we have
$$\sum_{k=1}^{m_n-d}2^{k}r^{-p(k-1)}n^{(p/2-\alpha)\vee (1+\nu-p/2)}\leqslant \left(2r^{-p}\right)^{m_n}n^{p/2-\alpha}\leqslant 2^{(1+\nu)m_n}n^{p/2-\alpha}$$
\begin{equation}
\leqslant 2^{1+\nu}\,n^{1/2+3\nu/2+\nu^2+p/2-\alpha},\;\;\;n\in\mathbb{N},\;\;\;\nu\in(0,\;p/2-\alpha-(1-p/2)). \label{a14}
\end{equation}
To obtain the last inequality, we used the estimate $2^{m_n}\leqslant 2\,n^{1/2+\nu}$. Since $\alpha>(p+1)/2$, the right-hand side of (\ref{a14}) tends to zero as $n\rightarrow\infty$, provided $\nu>0$ is small enough. \par

Furthermore, the inequality $p>2\alpha/(\alpha-1)$ yields $\beta=(1-1/\alpha)p/2>1$. We have
$$\sum_{k=1}^{m_n-d}2^{d+k}\left(2^{-d-k}\right)^{(1-1/\alpha)p/2}r^{-p(k-1)}\leqslant 2^{d(1-\beta)}\sum_{k=1}^{\infty}(2^{1-\beta}r^{-p})^k\leqslant 2^{d(1-\beta)}\sum_{k=1}^{\infty}2^{(1-\beta+\nu)k}.$$
The latter expression converges to zero as $d\rightarrow\infty$ if $\nu>0$ is small enough. Thus, for any $\eta>0$ one can find such $d,N\in\mathbb{Z}$ that the expression (\ref{a10}) is less than $\eta$ for all $n>N$. This implies (\ref{a8}). $\square$ \par

\renewcommand{\refname}{\normalfont\selectfont

\end{document}